\documentclass[11pt,reqno]{amsart}
\usepackage{amsmath,amssymb,physics}
\usepackage{amsthm,amsopn,amssymb,a4wide}
\usepackage{enumitem}
\usepackage[mathscr]{eucal}
\usepackage{tikz}
\usepackage[colorlinks]{hyperref}
\usepackage[capitalize,nameinlink]{cleveref}
\hypersetup{urlcolor=blue, citecolor=red, linkcolor=blue}
\usepackage{color, bm, amscd, tikz-cd}
\usepackage{nicematrix} 
\newcommand{\ran}{\operatorname{ran}}

\newcommand{\C}{\mathbb{C}}
\newcommand{\Z}{\mathbb{Z}}
\newcommand{\R}{\mathbb{R}}
\newcommand{\N}{\mathbb{N}}
\newcommand{\D}{{\mathbb D}}

\newcommand{\B}{\mathscr{B}}
\renewcommand{\S}{\mathscr S}
\renewcommand{\H}{\mathcal{H}}
\newcommand{\K}{\mathcal{K}}
\newcommand{\E}{\mathcal{E}} 

\newcommand{\z}{\mathbf{z}}

\newcommand{\V}{\mathscr V}
\newcommand{\Tc}{\mathscr T}
\renewcommand{\t}{\textbf t}

\DeclareMathOperator{\slim}{s\!\cdot\! lim}

\newtheorem{thm}{Theorem}[section]

\newtheorem{corollary}[thm]{Corollary}
\newtheorem{lemma}[thm]{Lemma}

\newtheorem{proposition}[thm]{Proposition}

\newtheorem{definition}[thm]{Definition}
\newtheorem{remark}[thm]{Remark}

\newtheorem{example}[thm]{Example}

\numberwithin{equation}{section}

\def\textmatrix#1&#2\\#3&#4\\{\bigl({#1 \atop #3}\ {#2 \atop #4}\bigr)}
\def\dispmatrix#1&#2\\#3&#4\\{\left({#1 \atop #3}\ {#2 \atop #4}\right)}
\numberwithin{equation}{section}

\def\textmatrix#1&#2\\#3&#4\\{\bigl({#1 \atop #3}\ {#2 \atop #4}\bigr)}
\def\dispmatrix#1&#2\\#3&#4\\{\left({#1 \atop #3}\ {#2 \atop #4}\right)}
\usepackage{fancyhdr}
\begin{document}
	
	\title[]{On factorization of the shift semigroup}

	\author{Tirthankar Bhattacharyya}
	\address{Department of Mathematics, Indian Institute of Science, Bangalore 560012.}
	\email{tirtha@iisc.ac.in}
	
	\author{Shubham Rastogi}
	\address{Department of Mathematics, Indian Institute of Science, Bangalore 560012.}
    \email{shubhamrastogi4ever@gmail.com}
	
	\author{Kalyan B. Sinha}
	\address{Jawaharlal Nehru Centre for Advanced Scientific Research, Bangalore 560064.}
	\email{kbs@jncasr.ac.in}
	
	\author{Vijaya Kumar U}
	\address{Department of Mathematics,  SRM Institute of Science and Technology, Kattankulathur - 603203.} 
    \email{vijayaku@srmist.edu.in}

	\maketitle
	
	\renewcommand{\thefootnote}{\fnsymbol{footnote}}
	
	\footnotetext{MSC: Primary: 47D03, 47A65, 47A68.\\
		Keywords: The shift semigroup, B-C-L Theorem, Commutant, Factorization, Stolz region.}
	
	\begin{abstract}
		Let $\E$ be a finite dimensional Hilbert space. This note finds all factorizations of the right shift semigroup $\S^\E=(S_t^\E)_{t\ge 0}$ on $L^2(\R_+,\E)$ into the product of $n$ commuting contractive semigroups, i.e., characterizes all $n$-tuples of commuting semigroups $(\V_1,\V_2,...,\V_n)$ where $\V_i=(V_{i,t})_{t\ge 0}$ for $i=1,2,...,n$ are semigroups of contractions satisfying $V_{i,t}V_{j,t}=V_{j,t}V_{i,t}$  for all $i$ and $j$ and $S_t^\E=V_{1,t}V_{2,t}\cdots V_{n,t}$ for all $t\ge 0.$ The factorizations are characterized by tuples of self-adjoint operators $\underline{A}=(A_1,A_2,...,A_n)$ and tuples of positive contractions $\underline{B}=(B_1,B_2,...,B_n)$ on $\E$ satisfying certain conditions which are stated in \cref{thm:psi12}. One of the tools of our analysis is a convexity argument using the extreme points of the {\em Herglotz } class of functions 
		\[P:=\{f:\D\to \C  \text{ is analytic}, \Re{f}>0 \text{ and }f(0)=1 \}.\]
	\end{abstract}
	
	\section{Introduction}
	All Hilbert spaces in this note are over the complex field and are separable.
	Given a Hilbert space $\E$, let $H^2_\D(\E)$ be the {\em Hardy space} of $\E$-valued holomorphic functions on the open {\em unit disc} $\D$. Let $H^\infty_\D(\B(\E))$ be the Banach algebra of bounded holomorphic $\B(\E)$-valued functions on $\D$. A function $\pi$ in $H^\infty_\D(\B(\E))$ is called {\em inner} if the multiplication operator $M_\pi$ on $H^2_\D(\E)$ is an isometry. We shall call a tuple $(T_1, T_2, \ldots T_n)$ of bounded operators on a Hilbert space a {\em commuting factorization} of a bounded operator $T$ on the same Hilbert space if the $T_i$ commute and $T = T_1 T_2 \ldots T_n$. A famous theorem of Berger, Coburn and Lebow describes all commuting factorizations of the right shift operator $M_z$ on $H^2_\D(\E)$ into contractions.  We briefly describe this below.

	For an integer $n\ge 2$, a Hilbert space $\E$,  orthogonal projections $P_1,...,P_n$ on $\E$ and unitaries $U_1,...,U_n$ on $\E$, the triple  $\{\E, \underline{P}=(P_1,...,P_n), \underline{U}=(U_1,...,U_n)\},$ is called a {\em Berger-Coburn-Lebow (B-C-L) triple} if
	\begin{enumerate}
		\item $\underline{U}$ is a commuting factorization of $I_\E$,
		\item $P_j+U_j^*P_iU_j=P_i+U_i^*P_jU_i\le I_\E$ for $i\ne j;$ and 
		\item $P_1+U_1^*P_2U_1+U_1^*U_2^*P_3U_2U_1+\cdots +U_1^*U_2^*\cdots U_{n-1}^*P_nU_{n-1}\cdots U_2U_1=I_\E.$
	\end{enumerate}

	Given a B-C-L triple	$\{\E,\underline{P},\underline{U}\}$ the functions $\pi_j \in H^\infty_\D(\B(\E))$ defined by 
	\begin{equation}\label{defn:psi_j}
		\pi_j(z)=U_j(P_j^\perp+zP_j)\text{ for }j=1,2,...,n
	\end{equation}
	are inner. The tuple of isometries $(M_{\pi_1}, \ldots , M_{\pi_n})$ is a commuting factorization of $M_z^\E$ on $H^2_{\D}(\E)$	as can be seen directly by using the assumptions above.  We shall call the tuple $(M_{\pi_1},...,M_{\pi_n})$ a {\em model}. It is a remarkable fact these are all the commuting factorizations into contractions, see \cite{BCL, BDF}.
	
	If $\{\E, \underline{P}, \underline{U}\}$ and $\{\mathcal F, \underline{Q}, \underline{W}\}$ are any two B-C-L triples, then the corresponding models are unitarily equivalent if and only if there is a unitary $\Xi : \E \rightarrow \mathcal F$ intertwining $U_i$ with $W_i$ and $P_i$ with $Q_i$ for each $i$.  
	
	The model tuples also feature in an analogue of the Wold decomposition which states that a commuting tuple of isometries $(V_1, V_2, \ldots , V_n)$ on a Hilbert space $\mathcal H$ decomposes the space  uniquely into a direct sum $\H=\H_p\oplus \H_u$ satisfying
	\begin{enumerate}
		\item $\H_p$ and $\H_u$ are reducing for each $V_i$,
		\item $(V_1|_{\H_p},...,V_n|_{\H_p})$ is a model,
		\item    $(V_1|_{\H_u},...,V_n|_{\H_u})$ is a commuting tuple of unitaries.
	\end{enumerate}

	The restatement of the result above in terms of a semigroup considers the discrete semigroups  $(V_1^m)_{m=0}^\infty,...,(V_n^m)_{m=0}^\infty$. They are commuting in the sense that $V_i^m$ commutes with $V_j^m$ for every $i,j=1,2,...,n$ and $ m \ge 1$ and they are of the following form.
	\begin{equation*}
		V_j^m\simeq\begin{pNiceMatrix}[first-row,last-col,nullify-dots]
			\H_p & \H_u & \\
			M_{\pi_j^m} &0 & \H_p\\
			0& (V_j|_{\H_u})^m & \H_u
		\end{pNiceMatrix}
	\end{equation*}
	for $ m\in \Z_+$ and $ j=1,2,...,n.$ 
	
	In this note, we shall be concerned with continuous analogues of the factorization and of the structure. A {\em semigroup} will mean a strongly continuous $C_0$-semigroup of contractive linear operators on a Hilbert space indexed by the positive cone $\R^n_+$ of $\R^n$ for some $n\in\N.$ 	Given a semigroup $\V=(V_\t)_{\t\in \R^n_+}$ acting on a Hilbert space $\H$ and a subspace $\H_0$ of $\H$ which reduces $\V,$ i.e., $\H_0$ reduces $V_\t$ for all $\t\in \R^n_+,$ we denote by $\V|_{\H_0}$ the semigroup $(V_\t|_{\H_0})_{\t\in\R^n_+}.$ A semigroup $\V$ is said to be an {\em isometric/unitary semigroup} if $V_\t$ is an isometry/unitary operator on $\H$ for all $\t\in \R^n_+.$ It is said to be a {\em  completely non-unitary (c.n.u) semigroup} if there is no non-zero reducing subspace $\H_0$ of $\H$ for $\V,$ such that  $\V|_{\H_0}$ is a unitary semigroup.
	
	Just like the right shift operator $M_z$ on $H^2_\D(\E)$ is special among isometries (recall the remarkable Wold decomposition), there is a special isometric semigroup called the right shift semigroup.
	
	The {\em right shift semigroup} $\S^\E = (S_t^\E)_{t\ge 0}$ on $L^2(\R_+,\E)$ is defined by
	\[(S_t^\E f)x=\begin{cases}
		f(x-t) &\text{if } x\ge t,\\
		0 & \text{otherwise}
	\end{cases}\]
	for $t\ge 0.$ When $\E=\C$ we just write $\S= (S_t)_{t\ge 0}$ instead of $\S^\E.$
	
	A landmark result of Cooper, which is a far reaching generalization of Wold decomposition, shows that the right shift semigroup is indeed special.
	
	\begin{thm}[Cooper \cite{Cooper}]\label[thm]{Cooper}
		Let $\V=(V_t)_{t\ge 0}$ be a one-parameter $C_0$-semigroup of isometries on a Hilbert space $\H.$ Then  $\H$ uniquely decomposes as the direct sum
		\begin{equation}\label{CooperDecomp}
			\H=\H_p\oplus \H_u
		\end{equation}
		of reducing subspaces $\H_p$ and $\H_u$ of $\V,$ so that
		\begin{enumerate}
			\item $\V|_{\H_p}$ is a c.n.u isometric semigroup.
			\item $\V|_{\H_u}$ is a  unitary semigroup. 
		\end{enumerate}
		Moreover, any c.n.u one-parameter isometric semigroup is unitarily equivalent to the shift semigroup $\S^\E=(S_t^\E)_{t\ge 0}$ on $L^2(\R_+,\E)$ for some Hilbert space $\E.$ 
	\end{thm}	
	In fact, $\H_u=\bigcap_{t\ge 0}V_t(\H)$ in the above theorem. Cooper's theorem firmly establishes the special role of  $\S^\E$. A lucid proof of Cooper's theorem can be found in \cite{Sz. -Nagy}.
	
	
	Two one-parameter semigroups $\Tc_1=(T_{1,t})_{t\ge 0}$ and $\Tc_2=(T_{2,t})_{t\ge 0}$  are said to be {\em commuting} if
	\begin{equation*}
		T_{1,t}T_{2,t}=T_{2,t}T_{1,t} \text{ for all } t\ge 0.
	\end{equation*}
	We shall see in  \cref{InsteadOfShalit} that the above is equivalent to 
	$$T_{1,t}T_{2,s}=T_{2,s}T_{1,t} \text{ for all } t, s \ge 0.$$
	
	\begin{definition}Let $n\ge 2.$ Let $(\Tc_1,\Tc_2,...,\Tc_n)$ be an $n$-tuple of  commuting one-parameter semigroups of contractions on a Hilbert space $\H$. It is said to be an
		\begin{enumerate}
			\item  {\em n-isometric semigroup on $\H$} if each $\Tc_j$ is an isometric semigroup,
			\item {\em n-unitary semigroup on $\H$} if each $\Tc_j$ is a unitary semigroup.
		\end{enumerate}	 
	\end{definition}

	A tuple $(\Tc_1,\Tc_2,...,\Tc_n)$ of semigroups of contractions on $\H,$ is said to be a  {\em factorization} of a semigroup $\Tc$  (also $\Tc_1,\Tc_2,...,\Tc_n$ are said to {\em factorize} $\Tc$) if $\Tc_j$'s are commuting and
	\begin{equation}\label{defn:fact}
		T_t=T_{1,t}T_{2,t}\cdots T_{n,t}
		\text{ for all }t\ge 0.
	\end{equation}
	Note that if $\Tc$ is an isometric semigroup then each $\Tc_j$ must be an isometric semigroup.
	
	Let $\mathscr T=(T_t)_{t\ge 0}$ be a semigroup. A semigroup $\mathscr A=(A_t)_{t\ge 0}$ is said to be a {\em divisor} of the semigroup $\mathscr T$ if there is another semigroup $\mathscr B=(B_t)_{t\ge 0}$ such that $\mathscr A$ and $\mathscr B$ are commuting and 
	\begin{equation*}
		A_tB_t=T_t \text{ for all }t\ge 0.
	\end{equation*}

	{\em The principal aim of this note is to find all contractive factorizations of the right shift semigroup $(S_t^\E)_{t\ge 0}$.}

	The B-C-L type structure theorem for a tuple of commuting isometric semigroups is proved  in \cref{thm:BCL} in  \cref{B-C-L}. The route to factorization is through finding the contractive commutants of the right shift semigroup. To this end, we consider a unitary conjugation of the right shift semigroup on to a vector valued Hardy space, where it is a certain semigroup of multiplication operators. Then using a well-known result of Sarason (\cite{Sarason}), we describe the commutants in terms of certain multipliers.  This is done in \cref{thm:commutant}. {\em The results of  \cref{B-C-L,shift-semigp} hold when $\mathcal E$ is any separable Hilbert space.}
	
	Finally, when $\E$ is a finite dimensional Hilbert space, we completely describe the factorizations of the right shift semigroup on $L^2(\R_+,\E)$ in \cref{Factorization}.
	
	The authors are thankful to Prof. B. V. Rajarama Bhat for suggesting some of the questions and to Prof. Kaushal Verma for information on the Stolz regions.

	\section{The B-C-L theorem}\label[section]{B-C-L}
	The main result is \cref{thm:BCL} which is a far reaching generalization of Wold decomposition.
	
	The beginning of this section warrants a discussion of some background material all of which can be found in \cite[Chap III]{Nagy-Foias}. When $\E$ is one-dimensional, we denote $H^2_\D(\E)$ by $H^2$ and $H^\infty_\D(\B(\E))$ by  $H^\infty$. For any contraction $T$ on $\H$ and a function $u\in H^\infty$ with $u(z)=\sum_{n\ge 0}a_nz^n,$   the operators $u(rT):=\sum_{n\ge 0}a_nr^nT^n$  are defined for $0<r<1.$ 	For those functions $u\in H^\infty$ for which the strong limit $\slim_{ r \rightarrow 1^{-}}u(rT) $ exists, $u(T)$ is defined as
	\begin{equation*}
		u(T):=\slim_{ r \rightarrow 1^{-}}  u(rT).
	\end{equation*}
	The above association of the operator $u(T),$ to the function $u\in H^\infty$ is known as the {\em $H^\infty$ functional calculus} for $T$.
	
	If $A$ is the generator of a contractive semigroup $\Tc$, then $A-I$ is always boundedly invertible. The {\em cogenerator} $T$ of $\Tc$ is defined to be the bounded operator   $T=(A+I)(A-I)^{-1}.$  A contraction $T$ is the cogenerator for some semigroup of contractions if and only if $1$ is not an eigenvalue of $T.$ A semigroup is normal, self-adjoint, unitary or isometric if and only if its cogenerator is a normal, self-adjoint, unitary or isometric operator respectively.
	
	The semigroup of functions $\varphi_t\in H^\infty$ defined by 
\begin{align} \label{phit} \varphi_t(\lambda)=e^{t\frac{\lambda+1}{\lambda-1}} \text{ for }t\ge 0 \end{align} 
has a special role because if $T$ is the cogenerator of a semigroup $\Tc,$ then the semigroup $\Tc$ is given by (the $H^\infty$ functional calculus) $T_t=\varphi_t(T)$. Let $\sigma_p(T)$ denote the set of eigenvalues for a bounded operator $T$. If $T$ is a contraction such that $1\notin \sigma_p(T),$ then the $H^\infty$ functional calculus defines $\varphi_t(T)$  for all $t\ge 0$. Consequently, given a holomorphic $\B(\E)$-valued function $\psi$ with $\|\psi\|_{\infty} = \sup\{\|\psi(z)\| : |z| \le 1\} \leq 1$ and $1\notin \sigma_p(\psi(z))$ for any $z\in \D$, we have $\varphi_t(\psi(z))$ defined for all $t\geq 0$ and $z\in \D.$ In fact, we get that $\varphi_t\circ \psi \in H^\infty_\D(\B(\E))$ with $\|\varphi_t\circ \psi\|_{\infty}\leq1,$ for all $t\ge 0.$ We are now ready for the results.

	\begin{lemma}\label[lemma]{lem1}
		For $\psi\in H^\infty_\D(\B(\E)),$ with $\|\psi\|_{\infty}\leq1$, the following are equivalent:
		\begin{enumerate}
			\item $1\in \sigma_p(M_{\psi}).$
			\item $1\in \sigma_p(\psi(z))$  for some $z\in \D.$
			\item $1\in \sigma_p(\psi(z))$ for every $z\in \D.$
		\end{enumerate}
	\end{lemma}
	\begin{proof}
		$(1)\implies (2)$ is clear.
		
		$(2)\implies (3).$ Let $1\in \sigma_p(\psi(z_0))$ some $z_0\in \D,$ then there exists a non zero vector $x\in \E$ such that $\psi(z_0)x=x.$ Since $\|\psi\|_{\infty}\leq1,$ by a result in \cite{Brown-Douglas}, we have $\psi(z)x$ is a constant function of $z$. As $\psi(z_0)x=x,$ we have $\psi(z)x=x$ for all $z\in \D.$ Thus, $x$ is an eigenvector for all $\psi(z)$ with eigenvalue 1. (See also \cite{Brown-Douglas}).
		
		$(3)\implies (2)$ is trivial.
		
		$(2)\implies (1).$ Let $\psi(z_0)x=x$ for some non zero vector $x\in \E.$ By the above argument, we get $\psi(z)x=x$ for all $z\in \D.$ We are done by noting that $M_{\psi}(\mathbf{1}\otimes x)=\mathbf{1}\otimes x,$ where $\mathbf{1}\otimes x$ denotes the constant function in $H^2_\D (\E)$ with constant value $x.$
	\end{proof}
	
	It will be convenient to consider the class
	\begin{equation} \label{psi_property}
		\mathcal C = \{ \psi \in H^\infty_\D(\B(\E)) \text{ with } \|\psi\|_{\infty}\le 1 \text{ and } 1\notin \sigma_p(\psi(z)) \text{ for any } z\in \D\}.
	\end{equation}
	
	Let $\psi$ be from the class $\mathcal C$. Then by  \cref{lem1}, $1\notin \sigma_p(M_\psi).$ By \cite[Chap. III]{Nagy-Foias}, $M_\psi$ is the cogenerator of a contractive semigroup, which, as a semigroup of multiplication operators is as follows.
	\begin{lemma}\label[lemma]{lem2}
		Let $\psi$ be from the class $\mathcal C$ (\cref{psi_property}) and let $\varphi_t$ be as in \cref{phit}. Then,  $(M_{\varphi_t\circ \psi})_{t\geq 0}$ is the contractive semigroup on $H^2_\D(\E)$ whose cogenerator is $M_\psi.$
	\end{lemma}
	
	\begin{proof}
		This is a straightforward consequence of the $H^\infty$  functional calculus. By \cite[Chap. III]{Nagy-Foias}, the semigroup is $(\varphi_t(M_{\psi}))_{t\geq 0}$. Then, we have
		\begin{align*}
			\varphi_t(M_\psi)&=\slim_{ r \rightarrow 1^{-}}  e^{t(rM_\psi+I)(rM_\psi-I)^{-1}}\\
			&=\slim_{ r \rightarrow 1^{-}} e^{tM_{(r\psi+I)}{(M_{(r\psi-I)})}^{-1}}\\
			&=\slim_{ r \rightarrow 1^{-}} e^{M_{t(r\psi+I)(r\psi-I)^{-1}}} ~(\text{Note  that } (r\psi-I)^{-1}\in H^\infty_\D(\B(\E)).)  \\
			&=\slim_{ r \rightarrow 1^{-}} M_{e^{t(r\psi+I)(r\psi-I)^{-1}}} = M_{\varphi_t\circ \psi}.
		\end{align*}
		To see the last equality, let $f\in H^2_\D(\E)$ and $\lim_{r\rightarrow 1^{-}}M_{e^{t(r\psi+I)(r\psi-I)^{-1}}}f=g$ (say), i.e., $\lim_{r\rightarrow 1^{-}}e^{t(r\psi+I)(r\psi-I)^{-1}}f=g.$ Then for any $z\in \D$ and $x\in \E,$ we have $$\lim_{r\rightarrow 1^{-}}\langle e^{t(r\psi(z)+I)(r\psi(z)-I)^{-1}}f(z),x\rangle=\langle g(z),x\rangle.$$
		As $\varphi_t(\psi(z))= \slim_{ r \rightarrow 1^{-}} e^{t(r\psi(z)+I)(r\psi(z)-I)^{-1}},$ we have  $\varphi_t (\psi(z))(f(z))=g(z).$ Hence $ M_{\varphi_t \circ \psi}f=g.$ This completes the proof.
	\end{proof}
	
	\begin{lemma}\label[lemma]{lem3}
		Let $U,P\in \B(\E)$ be such that $U$ is a unitary and $P$ is a projection. Let $\pi(z)=U(P^\perp+zP)$ for $z\in \D.$ Then,
		$1\in \sigma_p(\pi(z))$ for some $z\in \D$ if and only if $Ux=x$ for some non-zero $x\in \ran(P^{\perp}).$
	\end{lemma}
	\begin{proof}
		The proof follows from \cref{lem1}.
	\end{proof}

\begin{definition}
Two semigroups $\Tc_1$ and $\Tc_2$ are said to be unitarily equivalent by the unitary $\Lambda$ if $T_{2,t}=\Lambda T_{1,t}\Lambda^*$ for $t\ge 0$.
\end{definition}
		
	The following lemma is a natural one. We are supplying the proof because we could not find the lemma in the literature.

	\begin{lemma}\label[lemma]{lem:cogen-sg}
		Suppose $\Tc_i$ is a contractive semigroup on $\H_i$ with the cogenerator $T_i $ for $i=1,2$ and $\Lambda:\H_1\to \H_2$ is a unitary. Then, $T_1$ and $T_2$ are unitarily equivalent by the unitary $\Lambda$ (i.e., $T_2=\Lambda T_1\Lambda^*$) if and only if $\Tc_1$ and $\Tc_2$ are unitarily equivalent by the unitary $\Lambda$.
	\end{lemma}
	
	\begin{proof}
		Let $\Lambda:\H_1\to \H_2$ be a unitary such that $T_2=\Lambda T_1\Lambda^*.$ We recall $\varphi_t$ from \cref{phit} and for any $t\ge 0$, define the family of functions $\{ \varphi_{t,r} : 0 \le r < 1\}$ on $\D$ by
		$$ \varphi_{t,r}(z) = e^{t\frac{rz+1}{rz-1}}.$$
		Thus, for any contraction $X$ and $0\le r<1,$  we have $\varphi_{t,r}(X)=e^{t(rX+I)(rX-I)^{-1}}.$   Then, for $t\ge 0$ we have
		\begin{align*}
			\Lambda T_{1,t} \Lambda^*&=\Lambda \varphi_t(T_1)\Lambda^*
			=\Lambda\left(\slim_{r\to 1^-}\varphi_{t,r}(T_1)\right)\Lambda^*\\
			&=\slim_{r\to 1^-}\Lambda \varphi_{t,r}(T_1)\Lambda^*
			=\slim_{r\to 1^-}\lim_{n\to \infty} \Lambda P_{n,r}(T_1)\Lambda^* \\
			&=\slim_{r\to 1^-}\lim_{n\to \infty} P_{n,r}(T_2)
			=\slim_{r\to 1^-}\varphi_{t,r}(T_2)
			=T_{2,t}.
		\end{align*}
		In the above, $(P_{n,r})_{n=1}^\infty$ is a sequence of polynomials converging uniformly in $\D$ to $\varphi_{t,r}$ for every $0 < r < 1$.
		
		Conversely, let $\Lambda$ be a unitary such that $T_{2,t}=\Lambda T_{1,t}\Lambda^*$ for all $t\ge 0.$ For $t> 0,$ let $\vartheta_t(\lambda)=\frac{\lambda-1+t}{\lambda-1-t}.$ Then
		\begin{align*}
			\vartheta_t(T_{2,t})&=(T_{2,t}-1+t)(T_{2,t}-1-t)^{-1}=(\Lambda T_{1,t}\Lambda^*-1+t)(\Lambda T_{1,t}\Lambda^*-1-t)^{-1}\\&=\Lambda(T_{1,t}-1+t)\Lambda^*\Lambda(T_{1,t}-1-t)^{-1}\Lambda^*\\&=\Lambda(T_{1,t}-1+t)(T_{1,t}-1-t)^{-1}\Lambda^*=\Lambda\vartheta_t(T_{1,t})\Lambda^*.
		\end{align*}
		Therefore, by \cite[Chap. III, Thm 8.1]{Nagy-Foias}, \begin{align*}
			T_2=\slim_{t\to 0+}\vartheta_t(T_{2,t})
			=\slim_{t\to 0+}\Lambda\vartheta_t(T_{1,t})\Lambda^*
			=\Lambda\slim_{t\to 0+}\vartheta_t(T_{1,t})\Lambda^*	
			=\Lambda T_1\Lambda^*.
		\end{align*}
		This completes the proof.
	\end{proof}
	\begin{remark} \label[remark]{InsteadOfShalit}
		The tools used in the proof above also give us that with the same notations of \cref{lem:cogen-sg}, but with $\H_1=\H_2$ (i.e., the semigroups act on the same space) the following are equivalent:
		\begin{enumerate}
			\item $T_{1,t}$ and $T_{2,t}$ commute for every $t \ge 0$.
			\item $T_1$ and $T_2$ commute.
			\item $T_{1,t}$ and $T_{2,s}$ commute for every $t,s \ge 0$.
		\end{enumerate}
		The equivalence of (2) and (3) also follows from the proof of  \cite[Theorem 2.1]{Shalit}.
	\end{remark}
	
	An $n$-tuple $(V_1,...,V_n)$ of isometries  on $\H$ is said to be  {\em  completely non-unitary (c.n.u)} if there is no non-zero jointly reducing subspace $\H_0$ of $\H$ for $(V_1,...,V_n),$ such that  $V_j|_{\H_0}$ is a unitary  for each $j=1,2,...,n.$
	
	A subspace $\H_0$ of $\H$ is said to be {\em jointly reducing} for an $n$-tuple of semigroups $(\V_1,...,\V_n)$ if $\H_0$ reduces all $\V_j$'s simultaneously.
	
	An $n$-isometric semigroup $(\V_1,...,\V_n)$  on $\H$ is said to be {\em  completely non-unitary (c.n.u)} if there is no non-zero jointly reducing subspace $\H_0$ of $\H$ for $(\V_1,...,\V_n),$ so that  $\V_j|_{\H_0}$ is a unitary semigroup for each $j=1,2,...,n.$

	Let $(\V_1,...,\V_n)$ be an $n$-isometric semigroup on $\H.$ Let $\mathfrak V$ denote the product semigroup $\mathfrak V=(V_{1,t}\cdots V_{n,t})_{t\ge 0}.$ 
	Suppose $V_j$ is the cogenerator of the semigroup $\V_j$ and $V$ is the cogenerator of the product semigroup $\mathfrak V.$ 
	Then  $(\V_1,...,\V_n)$ is c.n.u  if and only if  $(V_1,...,V_n)$ is c.n.u  if and only if the single isometry $V$ is c.n.u.

	Let $\{\E,\underline{P},\underline{U}\}$ be a B-C-L triple such that $U_j$ does not fix any non-zero vector from $\ran (P_j^{\perp})$ for $j=1,...,n$ where $\underline{P}=(P_1,...,P_n)$ and $\underline{U}=(U_1,...,U_n)$. Let $\pi_j\in H^\infty_\D(\B(\E))$ be as given in \cref{defn:psi_j} and let $\varphi_t$ be as in \cref{phit}. Then by \cref{lem3}, $1\notin \sigma_p(\pi_j(z))$ for any $z\in \D.$ By \cref{lem2}, $(M_{\varphi_t \circ \pi_j})_{t\geq 0}$ is the contractive semigroup whose cogenerator is $M_{\pi_j},$ for all $j=1,...,n.$
	Since $(M_{\pi_1},...,M_{\pi_n})$ is an $n$-isometry on $H^2_\D(\E),$ we have  $((M_{\varphi_t \circ \pi_1})_{t\geq 0},..., (M_{\varphi_t \circ \pi_n})_{t\geq 0})$ is an $n$-isometric semigroup on $H^2_\D(\E).$
	\begin{definition}
		If $\varphi_t$ is as in \cref{phit}, then the $n$-isometric semigroup
		$((M_{\varphi_t \circ \pi_1})_{t\geq 0},..., (M_{\varphi_t \circ \pi_n})_{t\geq 0})$ is called a {\em model $n$-isometric semigroup}.	
	\end{definition}
	\begin{remark}\label[remark]{unique-BCL}
		\begin{enumerate}
			\item Any model $n$-isometric semigroup is c.n.u.
			\item Two model $n$-isometric semigroups are unitarily equivalent if and only if their corresponding B-C-L triples are unitarily equivalent.
		\end{enumerate}
	\end{remark}
	The above B-C-L triple $\{\E,\underline{P},\underline{U}\}$ is referred as the B-C-L triple of the model $n$-isometric semigroup $((M_{\varphi_t \circ \pi_1})_{t\geq 0},..., (M_{\varphi_t \circ \pi_n})_{t\geq 0}).$
	
	Let $\K$ be a Hilbert space and $((W_{1,t})_{t \geq 0},...,(W_{n,t})_{t \geq 0})$  be an $n$-unitary  semigroup on $\K$. Then,  $((M_{\varphi_t \circ \pi_1}\oplus W_{1,t})_{t\geq 0},...,(M_{\varphi_t \circ \pi_n}\oplus W_{n,t})_{t\geq 0})$  is an $n$-isometric semigroup on $H^2_\D(\E)\oplus \K.$ The following  theorem  tells us that, up to a unitary equivalence, these are the only $n$-isometric semigroups.

	\begin{thm}\label[thm]{thm:BCL}
		Let $(\V_1,...,\V_n)$ be an $n$-isometric semigroup on a Hilbert space $\H.$ Then,  $\H$ decomposes uniquely as a direct sum $\H=\H_p\oplus \H_u$ of joint reducing subspaces $\H_p$ and $\H_u$ for $(\V_1,...,\V_n)$   such that
		\begin{enumerate}
			\item   $(\V_1|_{\H_p},...,\V_n|_{\H_p})$ is a c.n.u $n$-isometric semigroup.		
			\item    $(\V_1|_{\H_u},...,\V_n|_{\H_u})$ is an $n$-unitary semigroup.
		\end{enumerate}
		Moreover, any c.n.u $n$-isometric semigroup is unitarily equivalent to a model $n$-isometric semigroup.
	\end{thm}

	\begin{proof}
		Let $V_i$ denote the cogenerator of $\V_i$ for $i=1,...,n.$ By \cite[Chap. III]{Nagy-Foias} and \cref{InsteadOfShalit}, $(V_1,...,V_n)$ is an $n$-isometry. By the fundamental result of Berger, Coburn and Lebow, the space $\H$ decomposes uniquely into a direct sum of jointly reducing subspaces $\H=\H_p\oplus \H_u$ such that $(V_1|_{\H_p},...,V_n|_{\H_p})$ is a c.n.u $n$-isometry and    $(V_1|_{\H_u},...,V_n|_{\H_u})$ is an $n$-unitary. 
		
		As $\H_p$ reduces the cogenerator $V_j,$ it reduces the semigroup $\V_j.$ Moreover, the cogenerator of the semigroup $\V_j|_{\H_p}$ is $V_j|_{\H_p}.$ Since $(V_1|_{\H_p},...,V_n|_{\H_p})$ is a c.n.u $n$-isometry, $(\V_1|_{\H_p},...,\V_n|_{\H_p})$ is a c.n.u $n$-isometric semigroup. Similarly, $(\V_1|_{\H_u},...,\V_n|_{\H_u})$ is an $n$-unitary semigroup on $\H_u.$
		
		Let $(\V_1,...,\V_n)$ be a c.n.u $n$-isometric semigroup on $\H$ and let $V_i$ denote the cogenerator of $\V_i$ for $i=1,...,n.$ Then $(V_1,...,V_n)$ is a c.n.u $n$-isometry. Therefore, there exists a B-C-L triple  $\{\E,(P_1,...,P_n),(U_1,...,U_n)\}$ and a unitary isomorphism $\Delta:\H\rightarrow H^2_\D(\E) $ such that 
		\begin{equation*}
			\Delta V_j\Delta^*=M_{\pi_j}\text{ for } j=1,...,n
		\end{equation*} 
		where $\pi_j(z)=U_j(P_j^\perp+zP_j).$ Since $V_j$ is the cogenerator of $\V_j,$ $1\notin \sigma_p(V_j)$  \cite[Chap. III]{Nagy-Foias}. Hence $1\notin \sigma_p(M_{\pi_j})$ for  $j=1,...,n.$ By \cref{lem1} and \cref{lem3}, $U_j$ does not fix any non zero vector in  $\ran P_j^\perp.$ 	Using \cref{lem:cogen-sg} and \cref{lem2}, we get \begin{equation*}
			\Delta V_{j,t}\Delta^*=M_{\varphi_t \circ \pi_j}\text{ for }j=1,...,n.
		\end{equation*}  This shows that $(\V_1,...,\V_n)$  is unitarily equivalent to the model $n$-isometric semigroup $((M_{\varphi_t \circ \pi_1})_{t\ge 0},...,(M_{\varphi_t \circ \pi_n})_{t\ge 0}).$
		
		For uniqueness of the decomposition $\H=\H_p \oplus \H_u,$ note that $\H_p$ and $\H_u$ are respectively the shift and unitary parts of the Wold decomposition of the product $V_1V_2\cdots V_n.$ 
	\end{proof}
	The model $n$-isometric semigroup arising  in part (1) of \cref{thm:BCL}, is unique up to a unitary equivalence.  The B-C-L triple of the model $n$-isometric semigroup is also called {\em as the B-C-L triple of the $n$-isometric semigroup $(\V_1,...,\V_n).$}

	There is a connection between $n$-parameter semigroups and one-parameter semigroups. 
	Let $\Tc=(T_\t)_{\t\in\R^n_+}$ be a strongly continuous $n$-parameter semigroup of contractions on $\H.$ If we define $$T_{1,t}=T_{(t,0,...,0)},T_{2,t}=T_{(0,t,0,...,0)},...,T_{n,t}=T_{(0,...,0,t)}$$ for all $t\ge 0,$ then $\Tc_j=(V_{j,t})_{t\ge 0}$ is a one-parameter strongly continuous semigroup of contractions on $\H$ for $j=1,2,...,n.$ Moreover, the semigroups $\Tc_j$ commute. Conversely, given any $n$-tuple $(\Tc_1,...,\Tc_n)$ of commuting semigroups of contractions on $\H$, if we define $$T_\t=T_{1,t_1}T_{2,t_2}\cdots T_{n,t_n}\text{ for }\t=(t_1,t_2,...,t_n)\in \R^n_+$$ where  $\Tc_j=(T_{j,t})_{t\ge 0}$ for $j=1,2,...,n.$ Then $\Tc=(T_\t)_{\t\in \R^n_+}$ is a strongly continuous $n$-parameter semigroup of contractions on $\H.$ The following remark is not hard to prove.  
	\begin{remark}\label[remark]{c.n.u}
		Let $\V$ be an $n$-parameter semigroup of isometries on $\H.$ Let $(\V_1,\V_2,...,\V_n)$ be the corresponding $n$-isometric semigroup on $\H.$ Let $\mathfrak V$ denote the product of the  semigroups $\V_1,\V_2,...,\V_n,$ i.e., $\mathfrak V=(V_{(t,t,...,t)})_{t\ge 0}$ and let $V$ be its cogenerator. Then \\
		(1). \[\bigcap_{\t\in \R^n_+}V_\t(\H)=\bigcap_{t\in \R_+}V_{(t,t,...,t)}(\H)=\bigcap_{n\in\Z_+} V^n(\H).\] (2). $\V$ is c.n.u if and only if $(\V_1,\V_2,...,\V_n)$ is c.n.u if and only if $\mathfrak{V}$ is c.n.u if and only if $\bigcap_{\t\in \R^n_+}V_\t(\H)=\{0\}.$
	\end{remark}
	
	The c.n.u part of \cref{thm:BCL} reformulated in terms of $n$-parameter semigroup of isometries is as follows.
	\begin{corollary}
		Let $((M_{\varphi_t\circ \pi_1})_{t\ge 0},..., (M_{\varphi_t\circ \pi_n})_{t\ge 0})$ be a model $n$-isometric semigroup on $H^2_\D(\E).$ Define 
		\begin{equation}\label{n-model}
			M_\t=M_{\varphi_{t_1}\circ \pi_1}M_{\varphi_{t_2}\circ \pi_2}\cdots M_{\varphi_{t_n}\circ \pi_n}\text{ for }\t=(t_1,...,t_n)\in \R^n_+.\end{equation}
		Then $(M_\t)_{\t\in \R^n_+}$ is a c.n.u $n$-parameter semigroup of isometries on $H^2_\D(\E).$ Moreover, any c.n.u $n$-parameter semigroup  of isometries is unitarily equivalent to a semigroup $(M_\t)_{\t\in \R^n_+},$ where $M_\t$ is defined by \cref{n-model} for some model $n$-isometric semigroup $((M_{\varphi_t\circ \pi_1})_{t\ge 0},..., (M_{\varphi_t\circ \pi_n})_{t\ge 0}).$
	\end{corollary}

	We noticed the following in this section.
	\begin{remark}
		(1). The B-C-L triple $\{\E,(P_1,...,P_n),(U_1,...,U_n)\}$ of a c.n.u $n$-isometric semigroup is a complete unitary invariant. \\
		(2). For $n=2,$  $U_2=U_1^*$ and $P_2=I-U_1P_1U_1^*.$ \\
		(3). In contrast to the  c.n.u one-parameter semigroups of isometries where there are only countably many (indexed by the multiplicity of the shift semigroup), there are uncountably many non-isomorphic c.n.u $n$-parameter semigroups  of isometries (c.n.u $n$-isometric semigroups) if $n\ge 2.$ 
	\end{remark}
	
	\section{The right shift semigroup and its commutant}\label[section]{shift-semigp}
	
	The  cogenerator $S^\E$ of the right shift semigroup $\S^\E=(S_t^\E)_{t\ge 0}$ on $L^2(\R_+, \E)$ is a pure isometry. See \cite[p. 153]{Nagy-Foias} for an explicit expression of $S^\E$.
	Let  $W_\E:L^2(\R_+,\E)\to H^2_\D(\E)$ be  the unitary defined  by
	\begin{equation}\label{shift}
		W_\E({(S^\E)}^n(\sqrt{2}e^{-x}f))=z^nf\text{ for } f\in \E, n\ge 0.
	\end{equation}
	For this section and the next,  we need a notation: $\z^\E\in H^\infty_\D(\B(\E))$ will denote the operator valued function  $\z^\E(w)=wI_\E$  for $w\in \D$ where $I_\E$ denotes the identity operator on $\E.$ Then $W_\E S^\E W_\E^*=M_{\z^\E}$.  Hence from \cref{lem:cogen-sg} and \cref{lem2}, we see that
	\begin{equation}\label{conj}
		W_\E S_t^\E W_\E^*=M_{\varphi_t\circ \z^\E} \text{ for all } t\ge 0.
	\end{equation}

	Let $\psi$ be from the class $\mathcal C$ (\cref{psi_property}). Then  $(M_{\varphi_t\circ \psi})_{t\ge 0}$ is a contractive semigroup on $H^2_\D(\E)$ which commutes with $(M_{\varphi_t\circ \z^\E})_{t\ge 0}.$ The following theorem shows that, these are the only contractive semigroups which commute with $(M_{\varphi_t\circ \z^\E})_{t\ge 0}.$
	
	\begin{thm} \label[thm]{thm:commutant}
		Let $\varphi_t$ be as in \cref{phit}. Let $(T_t)_{t\ge 0}$ be a contractive semigroup on $H^2_\D(\E)$ which commutes with $(M_{\varphi_t\circ \z^\E})_{t\ge 0}.$ Then, there exists a unique  $\psi$ from the class $\mathcal C$ (\cref{psi_property}) such that $T_t=M_{\varphi_t\circ \psi}$ for all $t\ge 0.$
	\end{thm}
	\begin{proof}
		Let $T$ be the cogenerator of $(T_t)_{t\ge 0}.$  Then, $T$ commutes with $M_{\z^\E}$ as $M_{\z^\E}$ is the cogenerator of $(M_{\varphi_t\circ \z^\E})_{t\ge 0};$ see \cref{InsteadOfShalit}. Therefore, by \cite{Sarason},  $T=M_\psi$ for some function $\psi\in H^\infty_\D(\B(\E)).$  
		Since $T$ is the cogenerator of a contractive semigroup, we have that $\norm{\psi}_\infty\le 1$ and $1$ is not an eigenvalue of $T$; see \cite[Chap. III, Thm. 8.1]{Nagy-Foias}. Hence $1\notin \sigma_p(\psi(z))$ for any $z\in \D$ by \cref{lem1}. Finally, note that $T_t=M_{\varphi_t \circ \psi }$ for all $t\geq 0$ by \cref{lem2}.
		
		The uniqueness of $\psi$ follows from the fact that: For $\psi_1,\psi_2\in \mathcal C,$ the semigroups $(M_{\varphi_t\circ \psi_1})_{t\ge 0}$ and  $(M_{\varphi_t\circ \psi_2})_{t\ge 0}$ are equal iff the cogenerators $M_{\psi_1}$ and $M_{\psi_2}$ are equal which happens if and only if $\psi_1=\psi_2.$
	\end{proof}
	
	The Sz.-Nagy Foias theory tell us that if $(V_1,V_2,...,V_n)$ is an $n$-isometry on $\H,$ with one of them being pure, say $V_1$ is pure, then there is a unitary $\Gamma:\H\to H^2_\D(\ker (V_1^*))$ such that $\Gamma V_1\Gamma^*=M_{\z^{\ker (V_1^*)}}$ and $\Gamma V_j\Gamma^*= M_{\psi_j}$  for some commuting inner functions $\psi_j\in H^\infty_\D(\B(\ker (V_1^*))), j=2,3,...,n$. We get a similar result at the semigroup level, as a corollary of  \cref{thm:commutant}.

	\begin{corollary}
		Let $\varphi_t$ be as in \cref{phit}. Let $(\V_1,\V_2,...,\V_n)$ be an $n$-isometric semigroup on $\H,$ with one of them being c.n.u, say $\V_1$ is c.n.u. Then, there exists a Hilbert space $\E$ and a unitary $\Lambda: \H\to H^2_\D(\E)$ such that $\Lambda \V_1\Lambda ^*=(M_{\varphi_t\circ \z^\E})_{t\ge 0}$ and $\Lambda \V_j\Lambda^*=(M_{\varphi_t\circ \psi_j})_{t\ge 0},$  for some commuting inner functions $\psi_j,j=2,3,...n$  from the class $\mathcal C$ (\cref{psi_property}).
	\end{corollary}
	\begin{proof}
		This follows from \cref{Cooper}, \cref{conj} and \cref{thm:commutant}. Note that $\psi_j$ is an inner function, as $M_{\psi_j}$ is the cogenerator of the isometric semigroup $(M_{\varphi_t\circ \psi_j})_{t\ge 0}.$  The functions $\psi_j$'s are commuting because the semigroups $\V_j$'s are commuting.
	\end{proof}
	
	\section{Factorization of the right shift semigroup}\label[section]{Factorization}
	
	In this section and the next, the Hilbert space $\mathcal E$ is of finite dimension and we shall obtain all the factorizations of $\S^\E$ on $L^2(\R_+,\E)$. We focus our attention on the finite-dimensional case because only in this case, we have the complete result of finding all factorizations. In the infinite-dimensional case, the factorizations found in this section will work. However, we do not know whether those are all.	Let $(\V_1,\V_2,...,\V_n)$ be a factorization of $\S^\E$. Hence, $(\V_1,\V_2,...,\V_n)$ is $n$-isometric semigroup on $L^2(\R_+,\E)$  such that
	\begin{equation}\label{eq:V1tV2t=St}
		V_{1,t}V_{2,t}\cdots V_{n,t}=S_t^\E\text{ for all }t\ge 0.
	\end{equation}
	Clearly, $\V_j$ commutes with $\S^\E$ for $j=1,2,...,n.$ We shall now employ the unitary described in \cref{shift} to move from the space $L^2(\R_+,\E)$ to the space $H^2_\D(\E)$ so that we can use the results on commutants. Since $W_\E S_t^\E W_\E^*=M_{\varphi_t\circ \z^\E}$ for all $t\ge 0$, the contractive semigroup $(W_\E V_{j,t}W_\E^*)_{t\ge 0}$ commutes with $(M_{\varphi_t\circ \z^\E})_{t\ge 0}.$ Thus, by \cref{thm:commutant}, $ W_\E V_{j,t}W_\E^*=M_{\varphi_t\circ \psi_j}, t\ge 0$ for some  {\em inner functions}  $\psi_j\in H^\infty_\D(\B(\E)),$ with  $1\notin \sigma_p(\psi_j(z))$ for any $z\in \D.$ So \cref{eq:V1tV2t=St} is equivalent to \begin{equation*}
		M_{\varphi_t\circ \psi_1}M_{\varphi_t\circ \psi_2}\cdots M_{\varphi_t\circ \psi_n}=M_{\varphi_t\circ \z^\E} \text{ for all } t\ge 0
	\end{equation*}
where $\varphi_t$ is as in \cref{phit}.	

	For $\psi_j$ from the class $\mathcal C$ (\cref{psi_property}), we shall call the tuple $(\psi_1 , \psi_2,...,\psi_n)$ a {\em factorizing tuple} if the contractive semigroups $(M_{\varphi_t\circ \psi_1})_{t\ge 0},(M_{\varphi_t\circ \psi_2})_{t\ge 0},...,(M_{\varphi_t\circ \psi_n})_{t\ge 0}$ factorize $(M_{\varphi_t\circ \z^\E})_{t\ge 0}.$

	\begin{proposition}\label[proposition]{Prop:equiv}
		If $\psi_j$'s are from the class $\mathcal C$ (\cref{psi_property}) for $j=1,2,...,n,$ then $(\psi_1 , \psi_2,..., \psi_n)$ is a  factorizing tuple if and only if  $\psi_j(z)$'s are commuting for all fixed $z\in\D$ and 	\begin{equation}\label{eq:master}
			\sum_{j=1}^n(\psi_j(z)+I)(\psi_j(z)-I)^{-1}=(z+1)(z-1)^{-1}I
		\end{equation}
		for all $z\in \D.$ Moreover, in such a case,  $\psi_j$'s are inner functions.
	\end{proposition}
	\begin{proof}
		If $(\psi_1 , \psi_2,...,\psi_n)$ is a factorizing tuple, then  $(M_{\varphi_t\circ \psi_j})_{t\ge 0}$'s commute. Thus, the cogenerators  $M_{\psi_j}$'s  commute. This happens if and only if $\psi_j(z)$'s commute for every $z\in \D.$
		
		Now, for $\psi_j\in  \mathcal C,$ such that $\psi_j(z)$'s commute for all $z\in \D,$
		\begin{equation*}
			M_{\varphi_t\circ \psi_1}M_{\varphi_t\circ \psi_2}\cdots M_{\varphi_t\circ \psi_n}=M_{\varphi_t\circ \z^\E} \text{ for all } t\ge 0
		\end{equation*}
		if and only if
		\begin{equation*}
			(\varphi_t\circ \psi_1)(\varphi_t\circ \psi_2)\cdots(\varphi_t\circ \psi_n) =\varphi_t\circ \z^\E\text{ for all } t\ge 0
		\end{equation*}
		or equivalently,
		\begin{equation*}\label{eq:e}
			e^{t(\psi_1(z)+I)(\psi_1(z)-I)^{-1}}e^{t(\psi_2(z)+I)(\psi_2(z)-I)^{-1}}\cdots e^{t(\psi_n(z)+I)(\psi_n(z)-I)^{-1}}=e^{t(z+1)(z-1)^{-1}I}
		\end{equation*}for all  $z\in \D, t\ge 0,$
		if and only if 
		\begin{equation*}
			e^{t\sum_{j=1}^n(\psi_j(z)+I)(\psi_j(z)-I)^{-1}}=e^{t(z+1)(z-1)^{-1}I}\text{ for all } z\in \D, t\ge 0
		\end{equation*}
		if and only if
		\begin{equation*}
			\sum_{j=1}^n(\psi_1(z)+I)(\psi_1(z)-I)^{-1}=(z+1)(z-1)^{-1}I
		\end{equation*}
		for all $z\in \D.$
		
		The final statement is a consequence of the fact that any tuple of commuting contractions, whose product is an isometry, are necessarily isometries.
	\end{proof}
	For $n>2$, the commutativity of $\psi_j(z)$'s has to be assumed  in the backward direction of the theorem above while for $n=2,$ \cref{eq:master} implies that $\psi_1$ and $\psi_2$ are commuting (it is precisely showed in the proof of Lemma \ref{lemma:psi1}). The following example illustrates this. Let 
\begin{align}\label{zeta} \zeta(z)=\frac{z+1}{z-1} \text{ for } z\in \D.\end{align}
	\begin{example}
		Let $A_1$ and $A_2$ be two non-commuting bounded self-adjoint operators on a Hilbert space $\E.$ Let 
		$\psi_j(z)=(iA_j+I)(iA_j-I)^{-1}$ for $j=1,2$
		and $\psi_3(z)=(-i(A_1+A_2)+\zeta(z)+I)(-i(A_1+A_2)+\zeta(z)-I)^{-1} $ for $z\in \D.$  Then $\psi_1(z)$ does not commute with $\psi_2(z)$ for all $z\in \D.$ Note that $ \psi_j\in \mathcal{C}$ for all $j=1,2,3$ and 
		\begin{equation*}
			\sum_{j=1}^3(\psi_j(z)+I)(\psi_j(z)-I)^{-1}=\zeta(z)I  \text{ for all }z\in \D.
		\end{equation*}
		But the semigroups $(M_{\varphi_t\circ \psi_1})_{t\ge 0}$ and $(M_{\varphi_t\circ \psi_2})_{t\ge 0}$  do not commute where $\varphi_t$ is as in \cref{phit}. Hence $(\psi_1,\psi_2,\psi_3)$ is not a factorizing tuple.
	\end{example}	
	\begin{remark}\label[remark]{remark:factor} 	If $(\psi_1 ,\psi_2,...,\psi_n)$ is a factorizing tuple, then
		\begin{equation}\label{eq:psi1}
			0 \ge	\Re{(\psi_j(z)+I)(\psi_j(z)-I)^{-1}}\ge \Re{(z+1)(z-1)^{-1}I}
		\end{equation}
		for all $z\in \D$ and for $j=1,2,...,n$. This is because of \cref{eq:master} and the fact that $\Re{(A+I)(A-I)^{-1}}\le 0$ for any contraction $A\in \B(\E)$ which does not have $1$ as an eigenvalue.	
	\end{remark}	
	To find all the factorizing tuples, we shall first find all $\psi\in \mathcal{C}$ satisfying the inequality \eqref{eq:psi1}. We start with a family of examples:
	\begin{proposition}\label[proposition]{Thm:Converse}
		Consider a self-adjoint operator $A$ and a positive contraction $B$ on $\E$. Let $\zeta$ be as in \cref{zeta}. Then $(iA+\zeta(z)B-I)$ is invertible. If $\psi:\D\rightarrow\B(\E)$ is defined as
		$$\psi(z)=(iA+\zeta(z)B+I)(iA+\zeta(z)B-I)^{-1},$$
		then   $\psi$ is a member of $\mathcal C$ which  satisfies \eqref{eq:psi1}.
	\end{proposition}
	\begin{proof}
		Since $\Re(iA+\zeta(z)B)\leq 0,$ the spectrum of $iA+\zeta(z)B$ is contained in the left half plane. Therefore $(iA+\zeta(z)B-I)$ is invertible.
		
		Clearly, $\psi$ is analytic. For each $z\in \D,$ $\psi(z)$ is a contraction because
		\begin{align*}
			\|\psi(z)(iA+\zeta(z)B-I)(x)\|^2&=\|(iA+\zeta(z)B+I)(x)\|^2\\
			&=\|(iA+\zeta(z)B)x\|^2+\|x\|^2+2\Re(\zeta(z))\langle Bx,x\rangle\\
			&\leq\|(iA+\zeta(z)B)x\|^2+\|x\|^2-2\Re(\zeta(z))\langle Bx,x\rangle\\
			&=\|(iA+\zeta(z)B-I)(x)\|^2.
		\end{align*}
		Suppose $\psi(z)(iA+\zeta(z)B-I)(x)=(iA+\zeta(z)B-I)(x)$ for some $z\in \D$ and $x\in \E.$ This implies $(iA+\zeta(z)B+I)(x)=(iA+\zeta(z)B-I)(x).$ Hence $x=0.$ Therefore $1\notin \sigma_p(\psi(z)).$ So $\psi$ lies in class $\mathcal C.$
		
		Note that $(\psi(z)+I)=2(iA+\zeta(z)B)(iA+\zeta(z)B-I)^{-1}$ and $(\psi(z)-I)^{-1}=\frac{1}{2}(iA+\zeta(z)B-I).$ Therefore $(\psi(z)+I)(\psi(z)-I)^{-1}=iA+\zeta(z)B.$ Finally, since $B$ is a positive contraction and $\Re{\zeta(z)}\le 0,$  $\psi$ satisfies \eqref{eq:psi1}.
	\end{proof}	
	
	Let $\varphi_t$ be as in \cref{phit}. For any $\psi\in \mathcal{C},$ from \cref{thm:commutant} and \cref{remark:factor} it follows that, if the semigroup $(M_{\varphi_t\circ \psi})_{t\ge 0}$ is a divisor of $(M_{\varphi_t\circ \z^\E})_{t\ge 0}$ then $\psi$ satisfies the inequality \eqref{eq:psi1}.	The following lemma establishes the converse.
	\begin{lemma}\label[lemma]{lemma:psi1}
		Let $\psi_1$ be from the class $\mathcal C$ (\cref{psi_property}) and satisfy \eqref{eq:psi1}. Then, there exists a unique $\psi_2$ in $\mathcal C$ satisfying \eqref{eq:psi1} such that $(\psi_1 , \psi_2)$ is a factorizing pair.
	\end{lemma}
	
	\begin{proof}
		For $z\in \D,$ let
		\begin{equation}\label{eq:h(z)}
			h_1(z)=(z+1)(z-1)^{-1}I-(\psi_1(z)+I)(\psi_1(z)-I)^{-1}.
		\end{equation}
		The inequality \eqref{eq:psi1} implies that  $\Re{h_1(z)}\le 0.$  Hence $( h_1(z) -I )^{-1}$  is in $\B(\E)$ for all $z\in \D.$ Now define $\psi_2:\D\to \B(\E)$ by \begin{equation}\label{eq:psi2}
			\psi_2(z)=(h_1(z)+I)(h_1(z)-I)^{-1}.
		\end{equation}
		Clearly $\psi_2$ is analytic.
		
		Let $T\in \B(\E)$ be such that  $\Re{T}\le 0.$ Let $y\in \E.$ Then by the straight-forward computation one can see that $\norm{(T+I)y}\le \norm{(T-I)y}.$ We shall use this fact in the following computation.
		
		Fix a $z\in \D$. Then for any $x$ in $\E$, there is a $y$ in $\E$ such that $x=(h_1(z)-I)y$, where $h_1(z)$ is as given in \cref{eq:h(z)}. Using $\Re{h_1(z)}\le 0$, we get
		$$			\norm{\psi_2(z)x}=\norm{(h_1(z)+I)y} \le \norm{(h_1(z)-I)y} =\norm{x}.$$
		This implies that $\norm{\psi_2}_\infty\le 1.$
		
		Suppose $1$ is an eigenvalue of $\psi_2(z)$ for some $z\in \D.$ Then, there exists a $y\in \E$ such that  $x=(h_1(z)-I)y\neq 0$ and $\psi_2(z)x=x.$ This implies that
		$$(h_1(z)+I)y=\psi_2(z)x=(h_1(z)-I)y.$$
		Thus $y=0$ and hence $x=0.$ This is a contradiction. Therefore, $1\notin \sigma_p(\psi_2(z))$ for any $z\in\D.$
		
		By the definition of $\psi_2,$ it is clear that  $\psi_1$ commutes with $\psi_2.$ It is a straightforward computation to show that the pair $(\psi_1,\psi_2)$ satisfies \cref{eq:master}. The uniqueness of $\psi_2$ follows from \cref{eq:master}.
	\end{proof}

	
	Given a $\psi_1$ in $\mathcal C$ satisfying \eqref{eq:psi1}, the unique $\psi_2$ obtained in \cref{lemma:psi1} is called the {\em factorizing partner} of $\psi_1$. We take recourse to some classical convex analysis to find the structure of a $\psi$ in $\mathcal C$ which satisfies \eqref{eq:psi1}.	
	Recall from the abstract that \begin{equation}
		P:=\{f:\D\to \C \text{ is analytic}, \Re{f}>0 \text{ and }f(0)=1 \}.
	\end{equation}
	It is clear that the set  $P$ is convex. It is  also compact in the topology of uniform convergence on compact subsets of  $\D.$
	\begin{thm}[see \cite{Holland}]\label[thm]{thm:extreme}
		The set of all extreme points of $P$ is precisely the set of functions
		\begin{equation}
			\left\{ \D \ni z \mapsto \frac{1+e^{i\theta}z}{1-e^{i\theta}z}: 0\le \theta< 2\pi\right\}.
		\end{equation}
	\end{thm}
	
	This will be our principal tool.
	\begin{proposition}\label[proposition]{tool}
		Let $f_j$  be holomorphic complex valued functions on $\D$ such that $\Re{f_j}\le 0$  for $j=1,2,...,n.$  Let $\zeta$ be as in \cref{zeta}. Suppose
		\begin{equation} \label{auxilliary_functions}
			f_1(z)+f_2(z)+\cdots +f_n(z)=\zeta(z)\text{ for all }z\in \D.
		\end{equation} Then $f_j(z)=i\Im{f_j(0)}-\Re{f_j(0)}\zeta(z)$ for all $z\in \D$ and $j=1,2,...,n.$
	\end{proposition}
	\begin{proof} Let $c_j=-\Re{f_j(0)}$ for $j=1,2,...,n.$ Note that   $$c_j\ge 0\text{ and } c_1+c_2+\cdots +c_n=1\text{ as } \zeta(0)=-1.$$   
		
		Let
		$$\Lambda=\{j\in \{1,2,...,n\}:c_j>0\}.$$
		If $j\in \{1,2,...,n\}\setminus \Lambda,$ $c_j=0,$ hence  by the open mapping theorem, $f_j$ is the constant function $f_j(0)$ and hence the desired conclusion holds trivially for  $j\in\{1,2,...,n\}\setminus \Lambda .$
		
		For $j\in \Lambda$ we have $\Re{f_j(z)}<0$ for all $z\in \D$ by the open mapping theorem.
		Let
		\begin{equation*}
			g_j(z)=\frac{1}{c_j}(i\Im{f_j(0)}-f_j(z)) \text{ for }  j\in \Lambda.
		\end{equation*} Hence, for $j\in \Lambda$ we have $\Re{g_j(z)}>0$ for all $z\in \D$  and $g_j(0)=1.$ Since $\sum_{j=1 }^n\Im{f_j(0)}=0$ and $f_j(z)=i\Im{f_j(0)}$ for $j\in\{1,2,...,n\}\setminus \Lambda,$ by \cref{auxilliary_functions} we have for $\zeta$ as in \cref{zeta},
		\begin{equation*}
			\sum_{j\in \Lambda}c_jg_j(z)=
			-\zeta(z) \text{ for all } z\in \D.
		\end{equation*}
		Therefore, by \cref{thm:extreme}
		we have $g_j(z)=-\zeta(z)$ for $j\in \Lambda.$ Hence
		\begin{equation*}
			f_j(z)=i\Im{f_j(0)}-\Re{f_j(0)}\zeta(z)\text{ for }j\in \Lambda \text{ and }  z\in \D.
		\end{equation*}	
		This completes the proof.
	\end{proof}
	\cref{tool} is proved for any $n\ge 2,$ however we use this result only for $n=2$ in the following lemma, to obtain  the structure of the divisors of $(M_{\varphi_t\circ \z^\E})_{t\ge 0}.$ This is perhaps the most important result of this section, it allows us to find the factorizing tuples completely.
	
	\begin{lemma}\label[lemma]{lemma:main}
		Let	$\psi$ be in $\mathcal C$ (\cref{psi_property}) and satisfy \eqref{eq:psi1}. Then there is a self-adjoint operator $A$ and a positive contraction $B$ on $\E$ such that for any $z\in\D,$
		\begin{equation}\label{psiAB}
			(\psi(z)+I)(\psi(z)-I)^{-1}=iA+\zeta(z)B
		\end{equation}	and
		\begin{equation*}
			\psi(z)=(iA+\zeta(z)B+I)(iA+\zeta(z)B-I)^{-1}.
		\end{equation*}
where $\zeta$ is as in \cref{zeta}. Moreover, the pair $(A,B)$ is unique.
	\end{lemma}
	\begin{proof}
		Let $\psi_1:=\psi.$ Let $\psi_2$ be its factorizing partner. Both of them satisfy
		$$\Re{(\psi_j(z)+I)(\psi_j(z)-I)^{-1}}\le 0$$
		for all $z\in \D$ and  $j=1,2$ by  \cref{remark:factor}. 
		
		For any unit vector $u\in \E$ and $j=1,2$ define
		\begin{equation*}
			f_{u,j}(z):=\langle (\psi_j(z)+I)(\psi_j(z)-I)^{-1}u,u\rangle \text{ for }z\in \D.
		\end{equation*}
		Then $f_{u,j}$ is analytic on $\D$ and  $\Re{f_{u,j}(z)}\le 0$ for all $z\in \D.$ Note for any unit vector $u\in\E$ that
		\begin{equation*}
			f_{u,1}(z)+f_{u,2}(z)=\zeta(z)\text{ for all }z\in \D.
		\end{equation*}
		Therefore, by \cref{tool} we get that $f_{u,j}(z)=i\Im{f_{u,j}(0)}-\Re{f_{u,j}(0)}\zeta(z).$ This implies that $$(\psi_j(z)+I)(\psi_j(z)-I)^{-1}=iA_j+\zeta(z)B_j,$$ where $A_j=\Im{(\psi_j(0)+I)(\psi_j(0)-I)^{-1}}$ and  $B_j=-\Re{(\psi_j(0)+I)(\psi_j(0)-I)^{-1}}.$  Clearly $A_j$'s are self-adjoint and $B_j$'s are positive. From \cref{eq:master} it follows that $A_1+A_2=0$ and $B_j$'s are contractions such that $B_1+B_2=I.$	Now
		\begin{align}
			\psi_j(z)&=[(\psi_j(z)+I)(\psi_j(z)-I)^{-1}+I][(\psi_j(z)+I)(\psi_j(z)-I)^{-1}-I]^{-1}\nonumber\\
			&=(iA_j+\zeta(z)B_j+I)(iA_j+\zeta(z)B_j-I)^{-1}\label{eq:psi-j-1-2}
		\end{align}
		for $z\in \D$ and $j=1,2.$ This completes the proof  with $A:=A_1$ and $B:=B_1$ as $\psi=\psi_1.$  The uniqueness of the pair $( A, B)$ is straightforward. In fact, $A=\Im{(\psi(0)+I)(\psi(0)-I)^{-1}}$ and $B=-\Re{(\psi(0)+I)(\psi(0)-I)^{-1}}.$
	\end{proof}
	\begin{remark}
		A semigroup $\V$ is a divisor of the semigroup $(M_{\varphi_t\circ \z^\E})_{t\ge 0}$ if and only if $\V$ is of the form $M_{e^{t(iA+\zeta(z)B)}}$ for some self-adjoint operator $A$ and a positive contraction $B$ on $\E$ where $\zeta$ is as in \cref{zeta}.
	\end{remark}	
	We have reached the main theorem of this note.

	\begin{thm}\label[thm]{thm:psi12}
		Given a factorizing tuple $(\psi_1,\psi_2,...,\psi_n)$, there exists a unique pair $(\underline{A},\underline{B})$ where $\underline{A}=(A_1,A_2,...,A_n)$ is a tuple of  self-adjoint operators and $\underline{B}=(B_1,B_2,...,B_n)$ is a tuple of  positive contractions on $\E$ satisfying
		\begin{enumerate}
			\item $A_jA_k=A_kA_j, B_jB_k=B_kB_j$ for $1\le j,k \le n-1,$
			\item $A_jB_k+B_jA_k=A_kB_j+B_kA_j$  for $1\le j, k \le n-1,$  and 
			\item $A_1+\cdots +A_n=0, B_1+\cdots+B_n=I$
		\end{enumerate} such that
		\begin{equation}
			(\psi_j(z)+I)(\psi_j(z)-I)^{-1}=iA_j+\zeta(z)B_j; \;\; j=1,2,..., n,
		\end{equation}	
		where $\zeta$ is as in \cref{zeta}. Consequently, the functions are of the form as in \cref{eq:psi-j-1-2}.
	\end{thm}
	\begin{proof}
		Since $(\psi_1,\psi_2,...,\psi_n)$ is a factorizing tuple each $\psi_j$ satisfies \eqref{eq:psi1}. Therefore, by \cref{lemma:main} there exists self-adjoint operators $A_1,A_2,...,A_n$ and positive contractions $B_1,B_2,...,B_n$ such that 
		\begin{equation*}
			(\psi_j(z)+I)(\psi_j(z)-I)^{-1}=iA_j+\zeta(z)B_j \text{ for all } j=1,2,..., n.
		\end{equation*}	
		By \cref{Prop:equiv}, $\psi_j$'s satisfy \cref{eq:master}. Hence, we have $A_1+\cdots +A_n=0$ and $B_1+\cdots+B_n=I.$
		
		Finally, as $(\psi_1,\psi_2,...,\psi_n)$ is a factorizing tuple, by \cref{Prop:equiv}, $\psi_j(z)$'s are commuting for all $j$. This is equivalent to say that $(\psi_j(z)+I)(\psi_j(z)-I)^{-1}$'s are commuting for all $j.$ This implies the conclusions (1) and (2). 
	\end{proof}
	
	Note that for $n=2$, the factorizing pair is indexed by a single self-adjoint operator $A$ and a positive contraction $B$. 
	
	In converse to \cref{thm:psi12}, if $\underline{A}=(A_1,A_2,...,A_n)$ is a tuple of  self-adjoint operators and $\underline{B}=(B_1,B_2,...,B_n)$ is a tuple of  positive contractions on $\E$ satisfying  (1)-(3) in \cref{thm:psi12}, and if $\psi_j$'s are defined as 
	\begin{equation}\label{def:psij}
		\psi_j(z)=(iA_j+\zeta(z)B_j+I)(iA_j+\zeta(z)B_j-I)^{-1}
	\end{equation}
	for all $j=1,2,...,n,$ then $\psi_j(z)$'s are commuting for all fixed $z$ and each $\psi_j$ is in class $\mathcal{C}.$ Moreover, $\psi_j$'s satisfies \cref{eq:master}. Hence $(\psi_1,\psi_2,...,\psi_n)$ is a factorizing tuple.

		\begin{remark} Let $A_1,A_2,...,A_{n-1}$ be self-adjoint operators and $B_1,B_2,...,B_{n-1}$ be  positive contractions on $\E.$ Then there exist a self-adjoint operator $A_n$ and a positive contraction $B_n$ on $\E$ such that $(\psi_1,\psi_2,...,\psi_n)$ is a factorizing tuple if and only the operators $A_j$ and $B_j$  satisfy conditions (1) and (2) in \cref{thm:psi12}, and $B_1+B_2+\cdots+B_{n-1}\le I,$ where $\psi_j$ are defined as in \cref{def:psij}. 
	\end{remark}

	The above result stated in terms of the semigroups is as follows:
	\begin{corollary}\label[corollary]{Factor}Let  $(\V_1,\V_2,...,\V_n)$ be a factorization of  $(M_{\varphi_t\circ\z^\E})_{t\ge 0}.$ Then, there exists a unique pair $(\underline{A},\underline{B}),$ where $\underline{A}=(A_1,A_2,...,A_n)$ is a tuple of  self-adjoint operators and $\underline{B}=(B_1,B_2,...,B_n)$ is a tuple of  positive contractions on $\E$ satisfying (1)-(3) of \cref{thm:psi12}, such that \begin{equation}\label{eq:Vjt}
			V_{j,t}=M_{\eta_{j,t}} \text{ for all }t\ge 0, j=1,2,...,n,
		\end{equation} where $\eta_{j,t}\in H^\infty_\D(\B(\E))$ is defined by \begin{equation}\label{eq:phi-j}\eta_{j,t}(z)=e^{t(iA_j+\zeta(z)B_j)} \text{ for }t\ge 0, j=1,2,...,n,
		\end{equation} where $\zeta$ is as in \cref{zeta}.
		
		Conversely, suppose $A_1,...,A_n$ are self-adjoint operators and $B_1,...,B_n$ are positive contractions on $\E$ satisfying (1)-(3) of \cref{thm:psi12}, if we define $\eta_{j,t}$ as in \cref{eq:phi-j} then $\eta_{j,t}\in H^\infty_\D(\B(\E))$ for $t\ge 0,j=1,2,...,n$  and  $(\V_1,\V_2,...,\V_n)$ as defined in \cref{eq:Vjt} is a factorization of  $(M_{\varphi_t\circ\z^\E})_{t\ge 0}.$
	\end{corollary}
	The above corollary, when expressed in terms of an $n$-parameter semigroup of contractions, can be reformulated as follows: It actually shows that if the diagonal $(V_{(t,t,...,t)})_{t\ge 0}$ of an $n$-parameter contractive semigroup $(V_\t)_{\t \in \R^n_+}$ is the shift semigroup, then it is an $n$-parameter semigroup of isometries determined by $n$ self-adjoint operators and $n$ positive contractions satisfying (1)-(3) of \cref{thm:psi12}:
	\begin{corollary}
		Let $\V=(V_\t)_{\t\in \R^n_+}$ be an $n$-parameter contractive semigroup such that 
		\begin{equation}\label{n-fact}
			V_{(t,t,...,t)}=M_{\varphi_t\circ\z^\E}\text{ for all }t\ge 0
		\end{equation}
where $\varphi_t$ is as in \cref{phit}. Then, there exists a pair $(\underline{A},\underline{B}),$ where $\underline{A}=(A_1,A_2,...,A_n)$ is a tuple of  self-adjoint operators and $\underline{B}=(B_1,B_2,...,B_n)$ is a tuple of  positive contractions on $\E$ satisfying (1)-(3) of \cref{thm:psi12}, such that
		\begin{equation}\label{n-para}
			V_{(t_1,t_2,...,t_n)}=M_{\eta_{1,t_1}\eta_{2,t_2}\cdots \eta_{n,t_n}}
		\end{equation}
		where $\eta_{j,t}$'s are $H^\infty_\D(\B(\E))$ functions as in \cref{eq:phi-j} for $t\ge 0, j=1,2,...,n.$

		Conversely, suppose $A_1,...,A_n$ are self-adjoint operators and $B_1,...,B_n$ are positive contractions on $\E$ satisfying (1)-(3) of \cref{thm:psi12}, if we define $\eta_{j,t}$ as in \cref{eq:phi-j} then $\eta_{j,t}\in H^\infty_\D(\B(\E))$ for $t\ge 0,j=1,2,...,n.$ Moreover, if we define $V_{(t_1,t_2,...,t_n)}$ as in \cref{n-para} for all $(t_1,t_2,...,t_n)\in \R^n_+$ then
		$(V_\t)_{\t\in \R^n_+}$ is an $n$-parameter semigroup of isometries satisfying \cref{n-fact}.
	\end{corollary}

	When $n=2,$ \cref{Factor} reduces as follows:
	\begin{remark}
		Let  $(\V_1,\V_2)$ be a factorization of  $(M_{\varphi_t\circ\z^\E})_{t\ge 0}$ where $\varphi_t$ is as in \cref{phit}. Then there exist  self-adjoint operators $A$ and $B$ in $\B(\E)$ with $0\le B\le I$ such that  \begin{equation*}
			V_{j,t}=M_{\eta_{j,t}} \text{ for all }t\ge 0, j=1,2,
		\end{equation*} where $\eta_{j,t}\in H^\infty_\D(\B(\E))$ is defined by \begin{equation*}\eta_{1,t}(z)=e^{t(iA+\zeta(z)B)}\text{ and }\eta_{2,t}(z)=e^{t(-iA+\zeta(z)(I-B))},
		\end{equation*} where $\zeta$ is as in \cref{zeta}.
		
		Conversely, suppose $A,B\in \B(\E)$ are self-adjoint and $0\le B\le I,$ if we define $\eta_{j,t}$ as in \cref{eq:phi-j} then $\eta_{j,t}\in H^\infty_\D(\B(\E))$ for $j=1,2,t\ge 0$ and  $(\V_1,\V_2)$ as defined in \cref{eq:Vjt} is a factorization of  $(M_{\varphi_t\circ\z^\E})_{t\ge 0}.$
	\end{remark}
	As an application of the factorization (\cref{Factor}) of the shift semigroup, we get another model for the c.n.u $n$-isometric semigroups, under some finite dimensional assumption.
	\begin{corollary}	
		Let $(\V_1,V_2,...,\V_n)$ be an $n$-isometric semigroup on $\H$. Let $V$ denote the cogenerator of the product semigroup $\V=\V_1\V_2\cdots\V_n.$ If $\dim (\ker V^*)<\infty,$ then upto a unitary equivalence, $\H$ decomposes as a direct sum $\H=\H_p\oplus \H_u$ of joint reducing subspaces $\H_p$ and $\H_u$ for $(\V_1,V_2,...,\V_n)$   such that
		\begin{enumerate}
			\item There is a unique triple $(\E, \underline{A}, \underline{B})$, where $\E$ is a finite dimensional Hilbert space, $\underline{A}=(A_1,A_2,...,A_n)$ is a tuple of  self-adjoint operators and $\underline{B}=(B_1,B_2,...,B_n)$ is a tuple of  positive contractions on $\E$ satisfying (1)-(3) of \cref{thm:psi12}; such that 
			$\H_p=H^2_\D(\E)$ and \begin{equation*}
				V_{j,t}|_{\H_p}=M_{e^{t(iA_j+\zeta(z)B_j)}} \text{ for all }t\ge 0, j=1,2,...,n.
			\end{equation*}
			\item    $\V_1|_{\H_u},\V_2|_{\H_u},...,\V_n|_{\H_u}$ are commuting unitary semigroups.
		\end{enumerate}
	\end{corollary}
	\begin{proof}It follows from applying \cref{Cooper} to the product semigroup $\V,$ \cref{conj} and \cref{Factor}. Note that if $\H_u$ is the unitary part in the Cooper's decomposition of the product semigroup $\V,$ it reduces all the semigroups $\V_1,\V_2,...,\V_n.$
	\end{proof}

	
	\section{The case of multiplicity one}
	It is instructive to take a look at the results obtained above when the Hilbert space $\E = \C$ because the complex analytic nature appears in full bloom. We shall show in this section that the factorizations of the right shift semigroup $\S$ are only the trivial ones, namely, those $(\V_1,\V_2,...,\V_n)$ where $V_{j,t}=e^{i a_j t}S_{b_jt}$ for some $a_j,b_j\in\R$ with $a_1+a_2+\cdots +a_n=0$ and $0\le b_j\le 1$  such that  $b_1+b_2+\cdots +b_n=1.$

	For any $a\in \R$ and $0\le b\le1,$ the isometric semigroup $(e^{i a t}S_{b t})_{t\ge 0}$ commutes with $(S_t)_{t\ge 0}.$ Therefore, by \cref{thm:commutant}, there exists $\psi\in H^\infty$ with $\norm{\psi}_\infty\le 1$ and $\psi\not\equiv 1 $ such that $We^{i a t}S_{b t}W^*=M_{\varphi_t\circ \psi}$ for all $t\ge 0,$ where $W:=W_\C$ is as given in \cref{shift}.  We explicitly find out the $\psi$ in the following lemma.
	\begin{lemma}\label[lemma]{lem:Salphapsi}
		Let $a$ be a real number and let $ b \in [0,1]$. Then there is a unique $\psi$ in the unit ball of $H^\infty$ of the unit disc such that  $We^{i a t}S_{b t}W^*=M_{\varphi_t\circ \psi}$. Moreover, the function $\psi$ is either  an automorphism or a constant.
	\end{lemma}
	
	\begin{proof} The uniqueness is clear. For the existence, we shall find out the $\psi$ explicitly. Define 
		$$\psi(z)=\alpha\frac{z-\beta}{1-\bar{\beta}z},$$
		with $\alpha=\frac{1+b+ia}{1+b-ia}$ and $\beta=\frac{1-b+ia}{1+b+ia} .$ It is an easy check that $|\alpha|=1,$ and $|\beta|<1$ if $b >0$ and $\beta=1$ if $b=0$. Note that $$\psi(z)=(ia+b\zeta(z)+1)(ia+b\zeta(z)-1)^{-1}.$$  Hence $\psi\in \mathcal C$ by \cref{Thm:Converse}. Since 
$$ia+b\zeta(z)=\frac{\psi(z)+1}{\psi(z)-1}\text{ for all } z\in \D,$$ 
we have
		$$e^{i a t}\varphi_{b t}(z)=e^{i a t}e^{b t\zeta(z)}=e^{t (ia +b \zeta(z))}=e^{t\frac{\psi(z)+1}{\psi(z)-1}}=\varphi_t\circ \psi(z)$$
		for all $z\in \D$ and $t\ge 0.$ Therefore, with $\varphi_t$ as in \cref{phit},
		$$M_{\varphi_t\circ \psi}=M_{e^{i a t}\varphi_{b t}}=e^{i a t}M_{\varphi_{b t}}=e^{i a t}W S_{b t}W^*=W e^{i a t}S_{b t}W^*$$ for $t\ge 0.$ This completes the proof.
	\end{proof}
	
\begin{figure}[h]
		\begin{tikzpicture}	
			\draw[color=black] (0,0) circle [radius=2];
			\draw[black,fill=gray!30,dashed] (0.6,0) circle [radius=1.4];
			\filldraw[black] (0.6,0) circle (1.5pt);
			\draw (0.6,0) node[below]{$\frac{M}{M+1}$};
			\draw[thick,<->] (-3,0) -- (3,0)  node[anchor=north west] {x axis};
			\draw (2.2,0) node[below]{$1$};
			\draw[thick,<->] (0,-3) -- (0,3) node[anchor=south east] {y axis};		
		\end{tikzpicture}
		\caption{A Stolz type region $\{ z \in \mathbb D : |1 - z|^2 < M(1 - |z|^2) \}.$} \label[figure]{fig:M1}
	\end{figure}
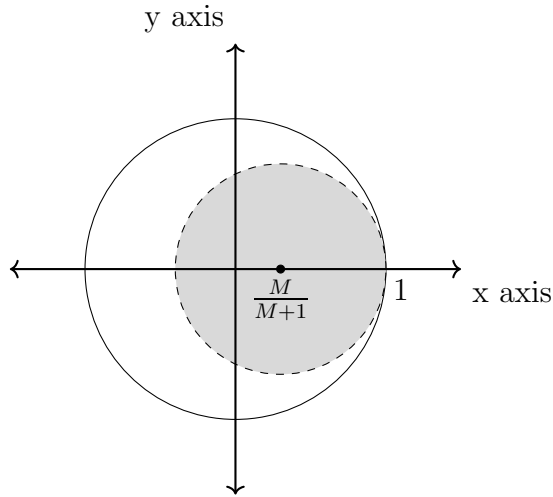

	The following complex analytic result is an easy exercise. It leads to a region (the shaded region in \cref{fig:M1}) which reminds an astute reader of the Stolz region stemming from the famous Stolz condition in Abel's limit theorem, see page 41 in \cite{Ahlfors}.
	
	\begin{thm}\label[thm]{thm:psiequiv}
		Let $\psi$ be a bounded holomorphic function with $\norm{\psi}_\infty\le 1$ and $\psi\not\equiv 1$. Then the following are equivalent
		\begin{enumerate}
			\item
			\begin{equation*}
				\Re{\frac{\psi(z)+1}{\psi(z)-1}}\ge \Re{\frac{z+1}{z-1}}\text{ for all } z\in \D.
			\end{equation*}
			\item \begin{equation*}
				\frac{| 1 - \psi(z)|^2}{1 - |\psi(z)|^2} \ge \frac{|1 - z|^2}{1 - |z|^2}.
			\end{equation*}
			\item $\psi$ is of the form $$\psi(z)=\alpha\frac{z-\beta}{1-\bar{\beta}z},\text{ for }z\in \D,$$ where $\alpha=\frac{1+b+ia}{1+b-ia}$ and $\beta=\frac{1-b+ia}{1+b+ia}$ for some $a$ and $b$ such that $a\in \R$ and $0\le b \le 1.$
			\item $\psi$ is  either a constant $e^{i \theta}$ for some $\theta \in (0,2\pi)$ or an automorphism such that $\psi^{-1}$ maps every Stolz type region $\{ z \in \mathbb D : |1 - z|^2 < M(1 - |z|^2) \}$ into itself.
		\end{enumerate}
	\end{thm}
	
	From the discussions above, there is a one-to-one correspondence between factorizing tuples $(\psi_1, \psi_2,...,\psi_n)$ and  \begin{equation}
		\left\{(\underline{a},\underline{b})\in \mathbb R^{n} \times \R_+^{n}:\sum_{j=1}^na_j=0\text{ and }\underline{b}\text{ is majorized by } (1,0,...,0)\right\}.
	\end{equation}(In $\mathbb R^n$, an $n$-tuple  $(x_1,...,x_n)$  is said to be majorized by another $n$-tuple $(y_1,...,y_n)$  if $\sum_{j=1}^kx_j\le \sum_{j=1}^ky_j$ for all $k=1,2,...,n$ and $\sum_{j=1}^nx_j= \sum_{j=1}^ny_j.$) Indeed, given such a pair $(\underline{a}=(a_1,...,a_{n}),\underline{b}=(b_1,...,b_{n})),$ if we define $$\alpha_j=\frac{1+b_j+ia_j}{1+b_j-ia_j},\quad \beta_j=\frac{1-b_j+ia_j}{1+b_j+ia_j}$$ for $j=1,2,...,n,$ and
	\begin{equation*}
		\psi_j(z)=\alpha_j\frac{z-\beta_j}{1-\bar{\beta_j}z},\text{ for }j=1,2,...,n
	\end{equation*}
	then $(\psi_1,\psi_2,...,\psi_n)$ is a factorizing tuple. As a consequence of \cref{remark:factor}, \cref{thm:psiequiv} and \cref{Prop:equiv}, it follows that every factorizing tuple is of this form. We close with the factorization promised at the beginning of the section.

	\begin{thm}
		If $(\V_1,\V_2,...,\V_n)$ is a factorization of $\S,$ then
		$$V_{j,t}=e^{i a_j t}S_{b_jt}$$ for a unique pair $(\underline{a},\underline{b})$ in $\mathbb R^{n} \times \R_+^{n}$ such that $\sum_{j=1}^na_j=0$  and $\underline{b}$ is majorized by $(1,0,...,0).$
	\end{thm}
	\begin{proof}
		We have $V_{j,t}=W^*M_{\varphi_t\circ \psi_j}W,$ for a factorizing tuple $(\psi_1, \psi_2,...,\psi_j),$ by \cref{thm:commutant} where $\varphi_t$ is as in \cref{phit}.  Every factorizing tuple is of the form
		$$\psi_j(z)=\alpha_j\frac{z-\beta_j}{1-\bar{\beta_j}z},\text{ for }j=1,2,...,n$$
		with 		$\alpha_j=\frac{1+b_j+ia_j}{1+b_j-ia_j}$ and $\beta_j=\frac{1-b_j+ia_j}{1+b_j+ia_j}$ for $j=1,2,...,n$ for a pair $(\underline{a},\underline{b})		$ as in the statement.  Now, by \cref{lem:Salphapsi}, we get that
		$$W^*M_{\varphi_t\circ \psi_j}W=e^{ia_jt}S_{b_j t}$$ for $j=1,2,...,n.$ This completes the proof.
	\end{proof}

	\noindent{\bf Acknowledgements}
	The first named author is supported by the J C Bose fellowship JCB/2021/000041 of SERB, the third named author is supported by the Senior Scientist Scheme of the Indian National Science Academy, and the foruth named author is suported by the  D S Kothari postdoctoral fellowship MA/20-21/0047. Research is supported by the DST FIST program - 2021 [TPN - 700661].\\

    \noindent{\bf Declarations}\\

	\noindent{\bf Data availability statement} Data sharing is not applicable to this article as
	no data sets were generated or analyzed during the current study.\\
	
	\noindent{\bf Conflict of interest} The authors declare that they have no conflict of interest.

\end{document}